\documentclass{pnastwo}

\usepackage{graphicx}

\usepackage{amssymb,amsfonts,amsmath}

\usepackage{mathptmx}
\renewcommand{\vec}[1]{\boldsymbol{#1}}

\contributor{Submitted to Proceedings of the National Academy of Sciences of the United States of America}
\url{www.pnas.org/cgi/doi/10.1073/pnas.0709640104}
\copyrightyear{2016}
\issuedate{Issue Date}
\volume{Volume}
\issuenumber{Issue Number}

\begin{document}


\title{The good, the bad, and the ugly: Bayesian model selection produces spurious posterior probabilities for phylogenetic trees} 






\author{
	Ziheng Yang
	\affil{1}{Department of Genetics, University College London, Gower Street, London WC1E 6BT, UK}
	\affil{2}{Radcliffe Institute for Advanced Studies, Harvard University, Cambridge, MA 02138, USA}
	\affil{3}{Institute of Applied Mathematics, Academy of Mathematics and Systems Science, Chinese Academy of Sciences, Beijing 100190, China}
	and
	Tianqi Zhu \affil{3}{Institute of Applied Mathematics, Academy of Mathematics and Systems Science, Chinese Academy of Sciences, Beijing 100190, China}
}

\contributor{Submitted to Proceedings of the National Academy of Sciences of the United States of America}

\significancetext{The Bayesian method is widely used to estimate species phylogenies using
	molecular sequence data.  While it has long been noted to produce spuriously high
	posterior probabilities for trees or clades, the precise reasons for this overconfidence
	are unknown.  Here we characterize the behavior of Bayesian model selection when the
	compared models are misspecified, and demonstrate that when the models are nearly equally
	wrong, the method exhibits unpleasant polarized behaviors, supporting one model with high
	confidence while rejecting others.   This provides an explanation for the empirical
	observation of spuriously high posterior probabilities in molecular phylogenetics.  
}


\maketitle 

\begin{article}


\begin{abstract}
	
The Bayesian method is noted to produce spuriously high posterior probabilities for
phylogenetic trees in analysis of large datasets, but the precise reasons for this
over-confidence are unknown.  In general, the performance of Bayesian selection of
misspecified models is poorly understood, even though this is of great scientific interest
since models are never true in real data analysis.  Here we characterize the asymptotic
behavior of Bayesian model selection and show that when the competing models are equally
wrong, Bayesian model selection exhibits surprising and polarized behaviors in large
datasets, supporting one model with full force while rejecting the others.  If one model
is slightly less wrong than the other, the less wrong model will eventually win when the
amount of data increases, but the method may become overconfident before it becomes
reliable.  We suggest that this extreme behavior may be a major factor for the spuriously
high posterior probabilities for evolutionary trees.
The philosophical implications of our results to the application of Bayesian model
selection to evaluate opposing scientific hypotheses are yet to be explored, as are the
behaviors of non-Bayesian methods in similar situations.

\end{abstract}









\section{Introduction}

\dropcap{T}he Bayesian method was introduced into molecular phylogenetics in the 1990s
\cite{Rannala1996, Mau1997, Li2000} and has since become one of the most popular methods
for statistical analysis in the field, in particular, for estimation of species
phylogenies \cite{Ronquist2012, Bouckaert2014, Lartillot2009, Chen2014}. It has been noted
that the method often produces very high posterior probabilities for trees or clades
(nodes in the tree).  In the first ever Bayesian phylogenetic calculation, a biologically
reasonable tree for five species of great apes was produced from a dataset of 11
mitochondrial tRNA genes (739 sites), but the posterior probability for that tree, at
0.9999, was uncomfortably high \cite{Rannala1996}. In the past two decades, the Bayesian
method has been used to analyze thousands of datasets, with the computation made possible
through Markov chain Monte Carlo (MCMC) \cite{Ronquist2012, Bouckaert2014}.  It has become
a common practice to report posterior clade probabilities only if they are $<100\%$
(because most estimates are 100\%). In some cases the high posterior probabilities are
decidedly spurious.  For example, conflicting trees may be inferred from the same data
under different evolutionary models.  Different trees may be
inferred depending on the species sampled in the dataset \cite{Bourlat2006}, or on whether
protein sequences or the encoding DNA sequences are analyzed \cite{Yang2008}.  In such
cases, the different trees cannot all be correct, even if the true tree is unknown.  The
concern is not so much that the inferred species relationships may be wrong as that they
are supported by extremely high posterior probabilities.

In the \emph{star-tree paradox}, large datasets were simulated using the star
tree, and then analyzed to calculate the posterior probabilities for the three binary
trees (Fig.~1).  Most biologists would like the posterior probabilities for
the binary trees to converge to $(\frac{1}{3}, \frac{1}{3}, \frac{1}{3})$ when the amount
of data increases \cite{Suzuki2002, Lewis2005, Yang2005}.  Instead they fluctuate among
datasets according to a statistical distribution, sometimes producing strong support for a
binary tree even though the data do not contain any information either for or against any
binary tree \cite{Steel2007, Yang2007, Susko2008}.

\begin{figure} [t]  \label{fig1trees}
	\centerline{\includegraphics[width=0.85\linewidth]{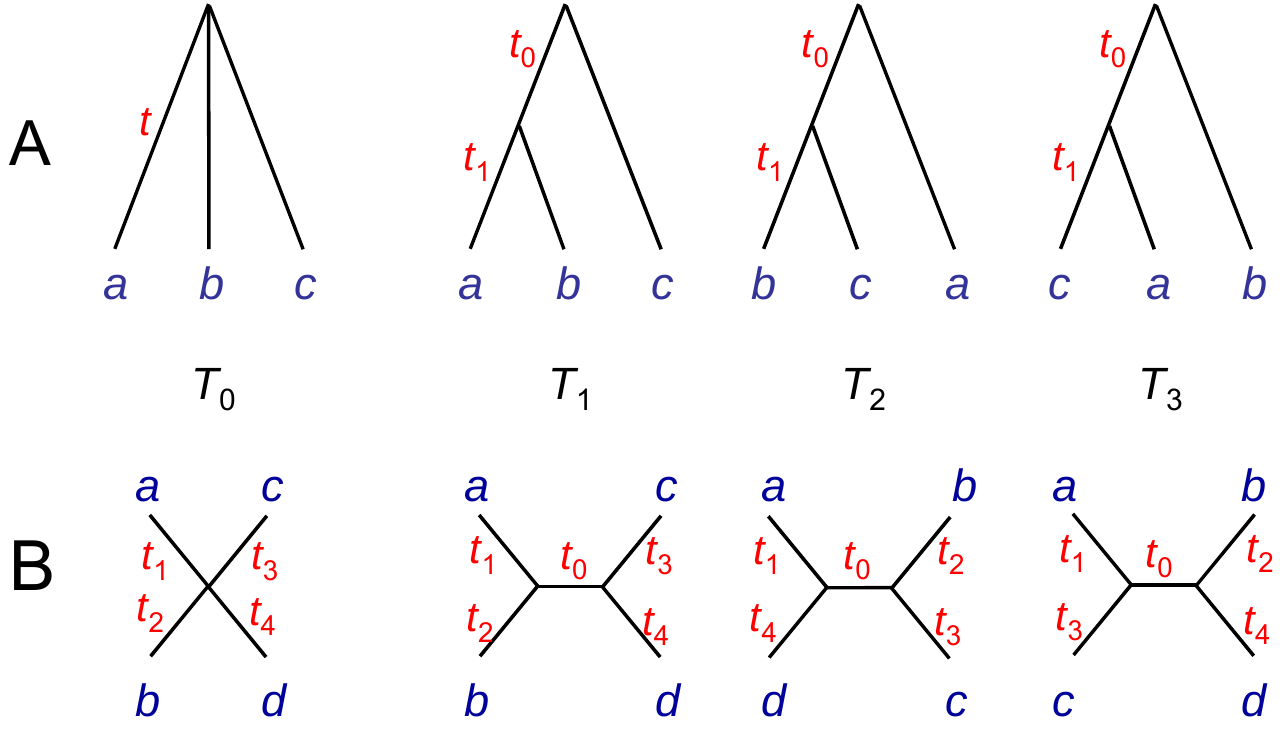}} 
	
	\caption{ (A) The three binary rooted trees for three species $T_1, T_2$, and $T_3$ and
		the star tree $T_0$. (B) The three binary unrooted trees for four species $T_1, T_2$, and
		$T_3$ and the star tree $T_0$. The branch length parameters are shown next to the
		branches, measured by the expected number of nucleotide changes per site.  In the
		star-tree simulations, the star tree is used to generate data, which are analyzed to
		calculate the posterior probabilities for the three binary trees, with the star tree
		excluded.} 
\end{figure}

Bayesian model selection is known to be consistent \cite{Dawid2011}. When the data size $n
\rightarrow \infty$, the true model `dominates', with its posterior probability
approaching 1.  If several models are equally right, the model with fewer parameters
dominates.  However, this theory applies only if the true model is included in the
comparison.  Given that a model is a simplified representation of the physical world, the
more common situation in real data analysis should be the comparison of models that are
all wrong.  Not many theoretical results appear to exist concerning Bayesian comparison of
misspecified models \cite{Berk1966}.

Here we study the asymptotic behavior of Bayesian model selection in a general setting
where multiple misspecified models are compared.  We are interested in how the posterior
probabilities for models behave when the data size increases.  Does the dynamics depend on
whether there are any free parameters in the models?  If one model is less wrong than
another (in a certain sense appropriately defined), will the less wrong model always win?
We present the proofs and mathematical analyses in SI Appendix.  In the main paper,
we summarize our results and illustrate them using three canonical simple problems.  Our
analysis suggests that the problem exposed by the star tree-paradox is actually far more
troubling than discussed previously \cite{Lewis2005, Yang2005, Steel2007, Yang2007,
	Susko2008}.

\section{Results}

\subsection{Problem description}

We consider independent and identically distributed (i.i.d.) models only.  The data
$\vec{x} = \{x_1, ..., x_n\}$ are an i.i.d.\ sample from the true model $g(\cdot)$.  We
consider two models as the case for more models is obvious.
Model $H_k$ has density $f_k(x|\theta_k)$, with $d_k$ free parameters ($\theta_k$), $k =
1, 2$.  We are in particular interested in models of the same dimension, with $d_1 = d_2 =
d$. In the Bayesian analysis, we assign a uniform prior for the two models ($\pi_1 = \pi_2
= \frac{1}{2})$ and also a prior for the parameters within each model $H_k$:
$f_k(\theta_k)$. The posterior model probabilities, $P_k = \mathbb{P}(H_k | \vec{x})$, are
then proportional to the marginal likelihoods: $M_k = f_k(\vec{x}) = \int f_k(\theta_k)
f_k(\vec{x}|\theta_k)\, \text{d}\theta_k$; that is, $P_1/P_2 = (\pi_1 M_1)/(\pi_2 M_2) =
M_1/M_2$. We are interested in the asymptotic behavior of $P_1$ in large datasets (as $n
\rightarrow \infty$).

The dynamics depends on how well the models fit the data.  Let $\hat{\theta_k}$ be the
maximum likelihood estimate (MLE) of $\theta_k$ under model $H_k$ from dataset $\vec{x}$. 
Let $\theta_k^*$ be the limiting value of $\hat{\theta_k}$ when the data size $n
\rightarrow \infty$.  In other words $\theta_k^*$ minimizes the Kullback-Leibler (K-L)
divergence from model $H_k$ to the true model
\begin{equation} \label{eq:KL}
  D_k = D_{\text{KL}}(g, f_k) = \int g(x) \log g(x) \,\text{d}x - \int g(x) \log f_k(x|\theta_k^*) \,\text{d}x,
\end{equation}
and is known as the \emph{best-fitting} or \emph{pseudo-true} parameter value under the
model \cite{White1982}.  $D_k$ (calculated at $\theta_k^*$) measures the distance from
$H_k$ to the true model, with $D_k \ge 0$.  We say a model is `right' if it encompasses
the true model, with $D = 0$, and `wrong' if $D > 0$.  Model 1 is less wrong than model 2
if $D_1 < D_2$.  Both models are `equally right' if $D_1 = D_2 = 0$ and `equally wrong' if
$D_1 = D_2 > 0$.

\subsection{Characterization of Bayesian model selection}

The asymptotic behavior of $P_1 = \mathbb{P}(H_1 | \vec{x})$ when $n \rightarrow \infty$
is analyzed in SI Appendix and summarized in Fig.~2.  We identify
three types of asymptotic behaviors: type-1 (`good'), type-2 (`bad') and type-3
(`ugly'), as defined below.  We also refer to three types of inference problems that
give rise to those behaviors.

Type 1 (`good') is for the posterior model probability $P_1$ to converge (as $n
\rightarrow \infty$) to a single reasonable value that is different from 0 and 1, such as
$\frac{1}{2}$.  In other words, in essentially every large dataset, $P_1 \approx
\frac{1}{2}$.  This behavior occurs when the two models are essentially
identical.  Examples include comparison of two identical models with no parameters, such
as $H_1: p = 0.5$ and $H_2: p = 0.5$ irrespective of the true $p$ in a coin-tossing
experiment (Fig.~2 cases $A_1$ and $A_2$), and overlapping models where
the best-fitting parameter values lie in the region of overlap (Fig.~2,
$A_3$ and $A_4$).  Whether the two models are both right ($A_1$ and $A_3$) or
both wrong ($A_2$ and $A_4$) does not affect the dynamics.  The case of overlapping
models is interesting.  If the truth is $p = \frac{1}{2}$ while the two compared models
are $H_1: 0.4 < p < 0.6$ and $H_2: 0 < p < 1$, and if we assign a uniform prior on $p$ in
each model, then as $n \rightarrow \infty$, $P_1 \rightarrow \frac{1}{1+0.2} =
\frac{5}{6}$, which appears more reasonable than $\frac{1}{2}$ as it favors the
more-informative model $H_1$.  At any rate, the comparison of identical or overlapping
models is unusual for testing scientific hypotheses.  This type of problem is not
considered further.

Type 2 ('bad') is for $P_1$ to converge to a nondegenerate
statistical distribution, such as $U(0, 1)$.  In other words, if we
analyze different large datasets, all generated from the same true
model, to compare two equally right or equally wrong models, $P_1$
varies among datasets according to a nondegenerate distribution. This
behavior occurs when the two compared models become unidentifiable as
the data size $n \rightarrow \infty$.  There are two scenarios.  In
the first, both models are right, with $D_1 = D_2 = 0$ (Fig.~2, $B_1$
and $B_2$). In the second both models are equally wrong (with $D_1 =
D_2 > 0$) but indistinct (Fig.~2, $B_3$ and $B_4$). We say that two
models are indistinct if and only if they, each at the best-fitting
parameter values, are unidentifiable, with $f_1(x|\theta_1^*) =
f_2(x|\theta_2^*)$ for essentially all $x$.  In other words, in
infinite data, the two models make essentially the same predictions
about the data and are unidentifiable.  In both scenarios of equally
right and equally wrong models, $P_1$ varies among datasets according
to a non-degenerate distribution.

Type 3 (`ugly') is for $P_1$ to have a degenerate two-point distribution, at values 0
and 1.  If we analyze large datasets to compare two models, we favour model 1 with total
confidence in some datasets and model 2 with total confidence in others.  This behavior is
observed when the two models are equally wrong and also distinct.  

It is remarkable that the asymptotic behavior is determined by whether
or not the compared models are distinct, and not by whether they are
both right or both wrong, or by whether the compared models have
unknown parameters.  For example, cases $B_1$ (two right models) and
$B_3$ (two equally wrong models) in Fig.~2 show the same `bad'
behavior, while cases $C_1$ (no free parameters) and $C_2$ (with free
parameters) show the same `ugly' behavior.

\begin{figure} [t]  \label{fig2summary}
	\centerline{\includegraphics[width=0.9\linewidth]{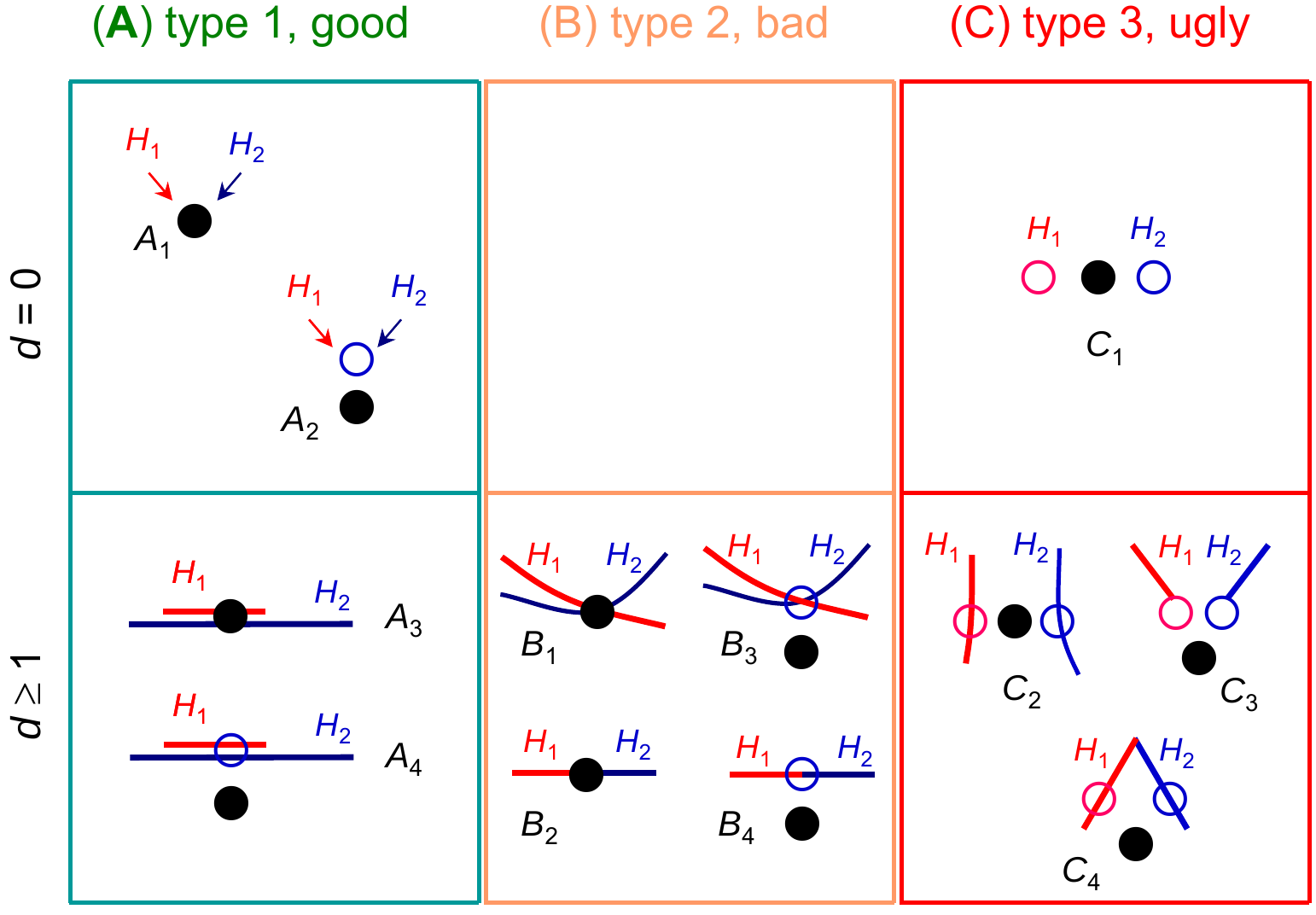}} 
	
	\caption{Classification of Bayesian model-selection problems involving two equally right
		or equally wrong models, each with $d$ free parameters.  Solid circles represent the true
		model, while the lines represent the parameter space of the compared models, with the
		empty circles to be the best-fitting parameter value $(\theta^*)$.  The two models are
		equally right (with $D_1 = D_2 = 0$) if the solid and empty circles coincide, and equally
		wrong (with $D_1 = D_2 > 0$) if they are separate.  The models are `indistinct' if the
		two empty cycles coincide (as in A and B), and are `distinct' if they are separate (as in
		C).  The green, orange, and red boxes indicate the three different asymptotic behaviors of Bayesian model
		selection when the data size $n \rightarrow \infty$.
	}
\end{figure}

\subsection{Problem 1. Fair-coin paradox (equally wrong models with no free parameter)}

Consider a coin-tossing experiment in which the coin is fair with the probability of heads
$p = \frac{1}{2}$.  We use the data of $x$ heads in $n$ tosses to compare two models:
$H_1: p = 0.4$ (tail bias) and $H_2: p = 0.6$ (head bias).  The two models are equally
wrong.  We assign a uniform prior for the two models ($\frac{1}{2}$ each), and calculate
the posterior model probability $P_1 = \mathbb{P}(H_1|x)$.  This is a type-3 problem
(Fig.~2, $C_1$).

As the models involve no free parameters, the likelihood $(L)$ and marginal likelihood
$(M)$ are the same, given by the binomial probability for data $x$.  The posterior odds is
the likelihood ratio
\begin{equation} 
  \frac{P_1}{1 - P_1} = \frac{M_1}{M_2} = \frac{0.4^x \cdot 0.6^{n - x}}{0.6^x \cdot 0.4^{n - x}}
   = \left( \frac{0.4}{0.6} \right)^{2x - n}.
\end{equation}
When $n$ is large, $P_1$ tends to be extreme (close to 0 or 1).  
Indeed $\alpha < P_1 < 1 - \alpha$ if and only if 
$|2x - n| < B = \frac{\log\{\alpha / (1 - \alpha)\}} {\log\{0.4/0.6\}}$.  
If $n$ is large, $2x-n$ is approximately $\mathbb{N}(0, n)$, so that
\begin{equation} \label{eq:FainCoinPP1}
  \mathbb{P}\{|2x - n| < B\} \approx 1 - 2\Phi(- \frac{B}{\sqrt{n}}) \approx \frac{2 B}{\sqrt{2\pi n}},
\end{equation}
where $\Phi$ is the cumulative distribution function (CDF) for $\mathbb{N}(0, 1)$.  If
$\alpha = 1\%$, we have $B$ = 11.33296, so that only 11 data outcomes will give $P_1$ in
the range (0.01, 0.99), with $x - \frac{n}{2}$ to be $-5, -4, \cdots, 5$.  For $n = 10^3,
10^4, 10^5, 10^6$, we have $\mathbb{P}\{0.01 < P_1 < 0.99\}$ = 0.280, 0.090, 0.0286, and
0.0090 using the normal approximation of eq.~\eqref{eq:FainCoinPP1}, or 0.272,
0.0876, 0.0277, and 0.0088 exactly by the binomial distribution.  Thus in large datasets,
moderate posterior probabilities will be rare, and either $H_1$ or $H_2$ will
be favored with posterior $>0.99$.  When $n \rightarrow \infty$, $P_1$ has a degenerate
two-point distribution, taking the values 0 and 1, each half of the times. This is the
type-3 `ugly' behavior.  Note that there is no information either for or against
either model in the data. Fig.~3A(i) shows the distribution of $P_1$
for $n = 10^3$.

Fig.~3A(ii) shows the comparison of $H_1: p = 0.42$ against $H_2: p =
0.6$ when the truth is $p = 0.5$.  Here $H_1$ is less wrong and will
eventually dominate.  However, in large and finite datasets, the more
wrong model $H_2$ can often receive high support. For example, for $n
= 10^3$, non-extreme posterior probabilities in the range $0.01 < P_1
< 0.99$ occur for only 13 data outcomes, with $x$ to be 504-516, and
in $14.8\%$ of datasets, $x$ is greater than those values so that $P_2
> 0.99$.  Indeed over the whole range $36 \le n \le 11611$, the more
wrong model $H_2$ is strongly favored too often, with $\mathbb{P}(P_2
>0.99) > 0.01$.  \emph{The method becomes overconfident before it
becomes reliable.} It may be noted that such strong support for the
more wrong model occurs only when the two models are opposing each
other.  It does not occur if both models are wrong in the same
direction: in the comparison of $H_1: p = 0.4$ and $H_2: p = 0.42$
when the truth is $p = 0.5$, the less wrong model $H_2$ dominates in
the posterior.

\subsection{Problem 2. Fair-balance paradox (equally right models or equally wrong and indistinct models)}

The true model is $\mathbb{N}(0, 1)$, and we compare two models $H_1: \mathbb{N}(\mu,
1/\tau)$, $\mu < 0$ and $H_2: \mathbb{N}(\mu, 1/\tau), \mu > 0$, with $\tau$ given. The
data may represent measurement errors observed on a fair balance while the models claim
that the balance has an unknown negative or positive bias.  The best-fitting parameter
value (the MLE when the data size $n \rightarrow \infty$) is $\mu^* = 0$ in each model,
when the two models become identical (indistinct).  Thus the two models are
equally right if $\tau = 1$  (Fig.~2, $B_2$), and are equally
wrong if $\tau = 1/9$ or 9 (Fig.~2, $B_4$).

We assign a uniform prior on the two models ($\frac{1}{2}$ each), and $\mu \sim
\mathbb{N}(0, 1/\xi)$ with $\xi$ fixed, truncated to the appropriate range under each
model.  The data $(x)$, an i.i.d.\ sample from $\mathbb{N}(0, 1)$, can be summarized as
the sample mean $\bar{x}$.  It can be shown that the posterior model probability $P_1 =
\mathbb{P}\{H_1 | x\}$ varies among datasets according to the following density
\begin{equation}  \label{eq:FareBalance}
  f(P_1) =\frac{\sqrt{\tau + \xi/n}}{\tau} \cdot \exp \left\{ \frac{[\Phi^{-1}(P_1)]^2}{2} 
  	\left[ 1 - \frac{1}{\tau} - \frac{\xi}{n\tau^2} \right]
  	\right\},
\end{equation}
where $\Phi^{-1}$ is the inverse CDF for $\mathbb{N}(0, 1)$ (SI Appendix). 

Fig.~3B shows the density of $P_1$ for different values of precision
$(\tau)$, with $n = 10^3$.  If $\tau = 1$, the two models are equally
right, and $f(P_1) \rightarrow 1$ when $n \rightarrow \infty$ so that $P_1$ behaves like a
$\mathbb{U}(0, 1)$ random number \cite{Lewis2005, Yang2005}.  If $\tau < 1$, the assumed
variance $(1/\tau)$ is larger than the true variance, so that the distribution has a mode
at $\frac{1}{2}$.  If $\tau > 1$, the assumed variance is too small, and $P_1$ has a
U-shaped distribution.  If one overstates the precision of the experiment, one tends to
over-interpret the data and generate extreme posterior model probabilities.  In all three
cases ($\tau <1, =1, >1$), $P_1$ has a non-degenerate distribution.

\begin{figure} [t]   \label{fig3problems123}
	\centerline{\includegraphics[width=1.0\linewidth]{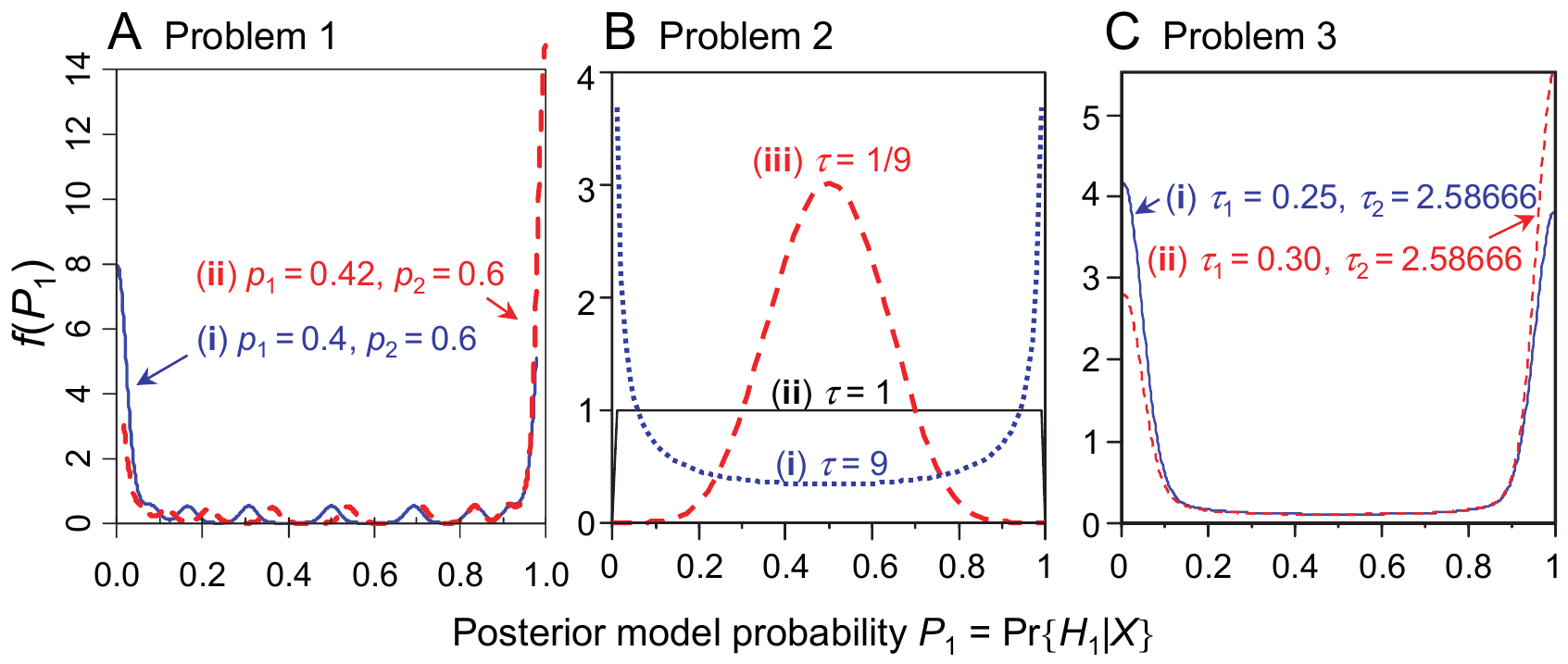}}
	
	\caption{The distribution of posterior model probability $P_1 = \mathbb{P}\{H_1|x\}$ in
		three inference problems. 	
		(A) Problem 1 (fair-coin paradox) is for a coin-tossing
		experiment, where the true model is $p = 0.5$ (a fair coin), and the compared models are
		(i) $H_1: p = 0.4$ and $H_2: p = 0.6$ so that the two models are equally wrong; and (ii)
		$H_1: p = 0.42$ and $H_2: p = 0.6$ so that $H_1$ is less wrong than $H_2$.  The data size
		(the number of coin tosses) is $10^3$. 
		(B) Problem 2 (fair-balance paradox) is for a normal-distribution example in which the
		true model is $\mathbb{N}(0, 1)$, and the two compared models are $H_1: \mathbb{N}(\mu,
		1/\tau), \mu < 0$ and $H_2: \mathbb{N}(\mu, 1/\tau), \mu > 0$, with variance $1/\tau$
		given.  The two models are equally right when $\tau = 1$ and equally wrong but indistinct
		when $\tau = 1/9$ or $9$.  The data size is $n = 10^3$. The plots for $n = 100$ or $\infty$
		are nearly the same. 
		(C) Problem 3 (fair-balance paradox) is for a normal-distribution
		example in which the true model is $\mathbb{N}(0, 1)$, and the two compared models are
		$H_1: \mathbb{N}(\mu, 1/\tau_1)$ and $H_2: \mathbb{N}(\mu, 1/\tau_2)$, with (i) $\tau_1 =
		0.25$ and $\tau_2 = 2.58666$, so that the two models are equally wrong and (ii) $\tau_1 =
		0.3$ and $\tau_2 = 2.58666$, so that $H_1$ is less wrong than $H_2$.  The prior is $\mu
		\sim \mathbb{N}(0, 1/\xi)$ under each model, with $\xi = 1$.  The data size is $n = 100$.
		All densities are estimated by simulating $10^5$ samples for $P_1$. }
\end{figure}

\subsection{Problem 3. Fair-balance paradox (equally wrong and distinct models)}

The true model is $\mathbb{N}(0, 1)$, and the two compared models are $H_1:
\mathbb{N}(\mu, 1/\tau_1)$ and $H_2: \mathbb{N}(\mu, 1/\tau_2)$, with $\tau_1 < 1 <
\tau_2$ given, while $\mu$ is a free parameter in each model.  The best-fitting parameter
value is $\mu^* = 0$ in each model, irrespective of the value of $\tau$ assumed.  Both
models are wrong because of the misspecified variance: $H_1$ is over-dispersed while
$H_2$ is under-dispersed.  They are equally wrong, in the sense that $D_1 = D_2$ in eq.~
\eqref{eq:KL}, if
\begin{equation}  \label{eq:BalanceKLcondition}
  \log \frac{\tau_1}{\tau_2} = \tau_1 - \tau_2
\end{equation}
(SI Appendix).  This is a type-3 problem (Fig.~2, $C_2$). We assign a
uniform prior over the models ($\frac{1}{2}$ each), and $\mu \sim
\mathbb{N}(0, 1/\xi)$, with $\xi$ given, within each model.  The
dataset, an i.i.d.\ sample of size $n$ from $\mathbb{N}(0, 1)$, can be
summarized as the sample mean $\bar{x}$ and sample variance $s^2 =
\frac{1}{n}\sum_i (x_i - \bar{x})^2$.  The posterior odds is given in
eq.~(15) in SI Appendix.

We use $\tau_1 = 0.25$ and $\tau_2 = 2.58666$, so that eq.~\eqref{eq:BalanceKLcondition}
holds and the two models are equally wrong, to generate independent variables $\bar{x}
\sim \mathbb{N}(0, 1/n)$ and $ns^2 \sim \chi^2_{n-2}$, and to calculate $P_1$.
Fig.~3C(i) shows the estimated density of $P_1$ for $n = 100$, with $\xi = 1$.  When $n \rightarrow \infty$,
$P_1$ degenerates into a 2-point distribution at 0 and 1, each with probability
$\frac{1}{2}$.  This is the same dynamics as in Problem 1
(Fig.~3A(i)), even though in Problem 1 the models do not involve any
unknown parameters while here they do.

Fig.~3C(ii) shows the density of $P_1$ when $\tau_1 = 0.3$ (which is closer to the true
$\tau = 1$ than is 0.25), so that $H_1$ is less wrong than $H_2$ (with $D_1 < D_2$).  In
this case when $n \rightarrow \infty$, $P_1 \rightarrow 1$.  However, in large but finite
datasets, $P_2$ for the more wrong model $H_2$ can be large in too many datasets: for
example, with $n = 100$, $\mathbb{P} \{P_2 > 0.99 \} = 0.0504$: in $5.04\%$ of datasets,
the more wrong model $H_2$ has posterior higher than $99\%$.

\subsection{Star-tree paradox and Bayesian phylogenetics}

In Bayesian phylogenetics \cite{Rannala1996, Mau1997}, each model has two components: the
phylogenetic tree describing the relationships among the species and the evolutionary
model describing sequence evolution along the branches on the tree
\cite{Yang1995YGF}.  Each tree $T_k$ has a set of time or branch length parameters
$(\vec{t}_k)$, which measure the amount of evolutionary changes along the branches.  The
evolutionary model may also involve unknown parameters $(\psi)$.  The tree and the
evolutionary model together specify the likelihood \cite{Felsenstein1981}, with $\theta =
\{\vec{t}, \psi\}$ to be the unknown parameters.  One of the trees is true, and all other
trees are wrong, while the evolutionary model may be misspecified.  The main objective is
to infer the true tree. The data consist of an alignment of sequences from the modern
species, and have a multinomial distribution in which the categories correspond to the
possible site patterns (configurations of nucleotides observed in the modern species)
while the data size is the number of sites or alignment columns \cite{Yang1994SB}.

Here we consider three simple cases involving 3 or 4 species (Fig.~1).  We
use the general theory described above to predict the asymptotic behavior of
posterior probabilities for trees and use computer simulation to verify the predictions.

Case A (Fig.~4A \& A') involves equally right models.  We use the rooted
star tree $T_0$ for three species with $t = 0.2$ (Fig.~1A) to generate
datasets to compare the three binary trees.  The Jukes-Cantor (JC) substitution model \cite{Jukes1969} is
used both to generate and to analyze the data, which assumes that the rate of change
between any two nucleotides is the same.  The molecular clock (rate constancy over time)
is assumed as well, so that the parameters in each binary tree are the two ages of nodes
($t_0, t_1$), measured by the expected number of nucleotide changes per site.

The best-fitting parameter values are $t_0^* = 0$ and $t_1^* = 0.2$ for each of the three
binary trees, in which case each binary tree converges to the true star tree.  We assign
uniform prior probabilities for the binary trees ($\frac{1}{3}$ each), and exponential
prior on branch lengths on each tree.  According to our characterization, this is a type-2
problem of comparing equally right models (Fig.~2, $B_2$), so the
posterior probabilities should have a nondegenerate distribution.  This case was
considered in previous studies \cite{Yang2005, Yang2007, Susko2008}, which generated
numerically the limiting distribution of the posterior probabilities for the binary trees
$(P_1, P_2, P_3)$ when $n \rightarrow \infty$, and pointed out that they do not converge
to $(\frac{1}{3}, \frac{1}{3}, \frac{1}{3})$ \cite{Lewis2005, Yang2005, Steel2007}.

Case B (Fig.~4B \& B') involves equally wrong models that are indistinct.
This is similar to case A except that the JC+$\Gamma$ model \cite{Jukes1969, Yang1993} is
used to generate data, with different sites in the sequence evolving at variable rates
according to the gamma distribution with shape parameter $\alpha = 1$.  The data are then
analyzed using JC (equivalently to JC+$\Gamma$ with $\alpha = \infty$).  The best-fitting
parameter values (i.e., the MLEs of branch lengths in infinite data) are $t_0^* = 0$ and
$t_1^* = 0.16441$ under each of the three binary trees.  The binary trees thus represent
equally wrong models (with $D_1 = D_2 = D_3 > 0$ in eq.~\eqref{eq:KL}) that are
indistinct.  The posterior tree probabilities have a non-degenerate distribution.  This is
the type-2 `bad' behavior for equal wrong and indistinct models (Fig.~2, $B_4$).

Case C (Fig.~4C \& C') involves equally wrong and distinct models. Like
case B, the simulation model is JC+$\Gamma$ with $\alpha = 1$, and the analysis model is
JC.  However, we do not assume the molecular clock and consider unrooted trees for four
species (Fig.~1B).  The true tree is the unrooted star tree $T_0$ of
Fig.~1B, with $t_1 = t_2 = t_3 = t_4 = 0.2$.  The best-fitting parameter
values (the MLEs of branch lengths in infinite data) are $t_0^* = 0.01037, t_i^* =
0.16409, i = 1, 2, 3, 4$, for each of the three binary trees (Fig.~1B).  As $t_0^* > 0$,
the three binary trees are different from the star tree and represent equally wrong and
distinct models (with $D_1 = D_2 = D_3 > 0$ in eq.~\eqref{eq:KL}). As this is a type-3
problem (Fig.~2, $C_4$), our theory predicts that as $n \rightarrow
\infty$, the posterior probabilities for the three binary trees should degenerate into a
three-point distribution, with probability $\frac{1}{3}$ each, for (1, 0, 0), (0, 1, 0),
and (0, 0, 1).  In other words, one of the binary trees will have posterior $\sim 100\%$
while the other two will have $\sim 0$. .  This is confirmed by simulation (table 1).

We note that most phylogenetic analyses involve unrooted trees as the clock assumption is
violated except for closely related species.  Furthermore, because of the violation of the
evolutionary model, all trees (or the joint tree-process models) represent wrong
statistical models.  Thus among the three cases considered in Fig.~4,
case C is the most relevant to analysis of real data, when Bayesian model selection
exhibits type-3 `ugly' behavior.  Previous analyses of the star-tree paradox
\cite{Yang2005, Yang2007, Susko2008} have deplored the `bad' behavior of Bayesian
phylogenetic method, but those studies examined case A only, so the real situation is
worse than previously realized.

A practically important scenario is where all binary trees are wrong because of violation
of the evolutionary model but the true tree is less wrong than the other trees.  We
present such a case in table 2, in which the data are simulated under JC+$\Gamma$ (with
$\alpha = 1$) using a binary tree with a short internal branch ($t_0 = 0.002$) and then
analyzed under JC. When the amount of data approaches infinity, the true tree will
eventually win, but there exists a twilight zone in which high posterior probabilities for
wrong trees occur too frequently; according to table 2, this zone is wider than $10^3 < n
< 10^5$. For example, at sequence length $n = 10^4$ and at the 1\% nominal level, the
error rate of rejecting the true tree is 25.0\% and the error rate of accepting a wrong
tree is 16.6\% (table 2).

\begin{figure} [t] \label{fig4startree}
	\centerline{\includegraphics[width=0.95\linewidth]{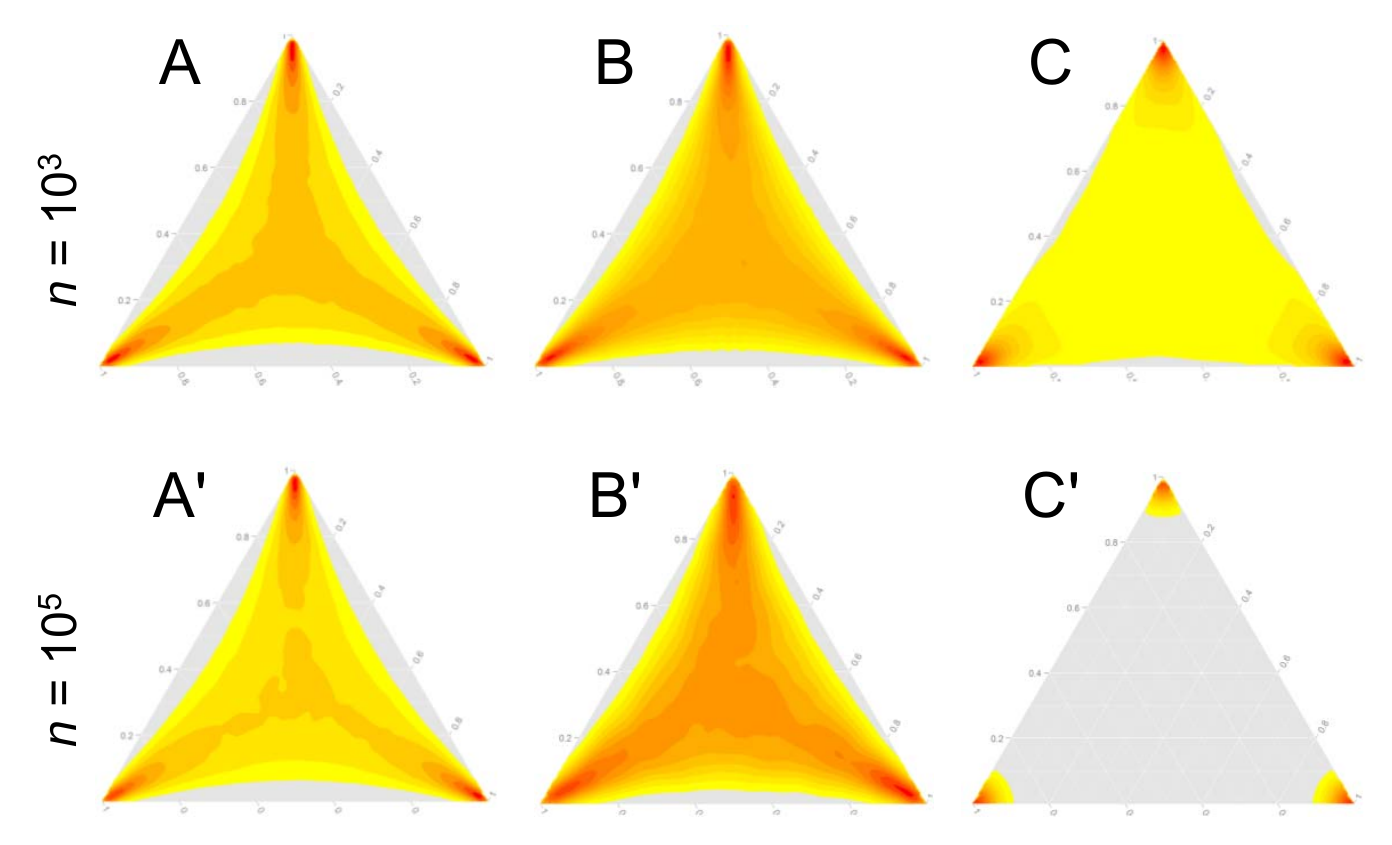}} 
	
	\caption{The distribution of posterior probabilities $(P_1, P_2, P_3)$ for the three
	binary trees $T_1, T_2$, and $T_3$ of Fig.~1, when datasets (sequence
	alignments of $n = 10^3$ or $10^5$ sites) are simulated using the star tree $T_0$ and
	analyzed to compare the three binary trees.  In (A) and (A'), the true tree is the star
	tree $T_0$ for three species of Fig.~1A, with $t = 0.2$.  Both the
	simulation and analysis models are JC, and the three binary trees are equally right
	models. In (B) and (B'), the true tree is the star tree $T_0$ for three species of
	Fig.~1A, with $t = 0.2$.  The simulation model is JC+$\Gamma$ (with $\alpha
	= 1$), and the analysis model is JC.  The three binary trees represent equally wrong and
	indistinct models. In (C) and (C'), the true tree is the star tree $T_0$ for four species
	of Fig.~1B, with $t_1 = t_2 = t_3 = t_4 = 0.2$.  The simulation model is
	JC+$\Gamma$ $(\alpha = 1)$ and the analysis model is JC.  The three binary trees
	represent equally wrong and distinct models.  The three corners in the plots correspond
	to points (1, 0, 0), (0, 1, 0), and (0, 0, 1), while the center is $(\frac{1}{3},
	\frac{1}{3}, \frac{1}{3})$. }
\end{figure}


\section{Discussion}

\subsection{High posterior probabilities for phylogenetic trees}

This work has been motivated by the phylogeny problem, and in particular by the empirical
observation of spuriously high posterior probabilities for phylogenetic trees
\cite{Suzuki2002, Lewis2005, Yang2005, Steel2007, Yang2007, Yang2008}.  We note that
certain biological processes such as deep coalescence \cite{Maddison1997, Xu2016}, gene
duplication followed by gene loss \cite{Nichols2001}, and horizontal gene transfer
\cite{Nichols2001, Maddison1997} may cause different genes or genomic regions to have
different histories.  However, as discussed in Introduction, posterior probabilities
for many trees or clades observed in real data analyses are decidedly spurious even if the
true tree is unknown.

One explanation for the spuriously high posterior probabilities for phylogenetic trees is
the failure of current evolutionary models  to accommodate interdependence among sites in
the sequence, leading to an exaggeration of the amount of information in the data.
Interacting sites may carry much less information than independent sites.  This
explanation predicts the problem to be more serious in coding genes than in noncoding
regions of the genome as noncoding sites may be evolving largely independently
due to lack of functional constraints.  However, empirical evidence points to the
opposite, with noncoding regions having higher substitution rates and higher information
content (if they are not saturated with substitutions), generating more extreme posteriors
for trees.  

Our results suggest that the problem may lie deeper and may be a
consequence of the polarized nature of Bayesian model selection when
all models under comparison are misspecified.  As the assumptions
about the process of sequence evolution are unrealistic, the
likelihood model is wrong whatever the tree, although the true tree
may be expected to be less wrong than the other trees.  As the
different trees constitute opposing models that are nearly equally
wrong, the inference problem is one of type-3 (Fig.~2, $C_4$).
Bayesian tree estimation may then be expected to produce extreme
posterior probabilities in large datasets.

\subsection{Bayesian selection of opposing misspecified models}

We have provided a characterization of model selection problems
according to the asymptotic behavior of the Bayesian method as the
data size $n \rightarrow \infty$ (Fig.~2 and SI Appendix).  While all
the problems considered here involve comparison of two `equally right'
or `equally wrong' models, three different asymptotic behaviors are
identified, which we label as type-1, type-2, and type-3.  The type-1
behavior is for the posterior model probability $P_1$ to converge to a
sensible point value, such as $\frac{1}{2}$.  We consider this to be a
`good' behavior, following phylogeneticists \cite{Suzuki2002,
Lewis2005, Yang2005}.  The rationale is that one would like a sure
answer given an infinite amount of data and the only reasonable sure
answer should be $\frac{1}{2}$ for each model, since the data contain
no information for or against either model.  This behavior occurs only
when the two models are identical or overlapping, a situation that
does not appear relevant to scientific inference.  With type-2
behavior, $P_1$ fluctuates among datasets (each of infinite size) like
a random number, so that strong support may be attached to a
particular model in some datasets.  Biologists were surprised at this
volatile behavior \cite{Suzuki2002, Lewis2005, Yang2005}.  
This occurs when the models are equally right or
equally wrong but indistinct.  In theory, type-2 behavior may not pose
a serious problem, because the parameter posteriors under the models,
if examined carefully, should make it clear that the competing models
essentially gave the same interpretation of the data and should lead
to the same scientific conclusion.  In data simulated in
\cite{Yang2005} or in Fig.~4A \& A', the estimates of $t_0$ should be
very close to 0, and all binary trees are similar to the same star
tree. Nevertheless this escaped our attention at the time.

With type-3 behavior, $P_1$ is $\sim 0$ in half of the datasets, and
$\sim 1$ in the other half.  This `ugly' polarized behavior occurs
when the two models are equally wrong and distinct.  Type-3 problems
may be the most relevant to practical data analysis given that all
models are simplified representations of reality and are thus wrong.
A variation to type-3 problems is when one model is only slightly less
wrong than another (Figs.~3A(ii), 3C(ii), table 2). While the less
wrong model eventually wins in the limit of infinite data, Bayesian
model selection is over-confident in large but finite datasets,
supporting the more wrong model with high posterior too often.

Note that the question of how the posterior model probability should
behave when large datasets are used to compare two equally wrong
models is somewhat philosophical and may not have a simple answer. One
position is to accept whatever behavior the Bayesian method
exhibits. This may be legitimate given that Bayesian theory is the
correct probability framework for summarizing evidence in the prior
and likelihood.  The polarized behavior in type-3 problems may then be
seen as a consequence of `user error' (for not including the true
model in the comparison), exacerbated by the large data size.  In this
regard we note that the posterior predictive distribution
\cite{Roberts1965, Box1980} can be used to assess the general adequacy
of any model or the compatibility between the prior and the
likelihood, and indeed this has been widely used to assess the
goodness of fit of models in phylogenetics \cite{Sullivan2005,
Rodrigue2006}.  Nevertheless, a number of sophisticated and
parameter-rich models have been developed for Bayesian phylogenetic
analysis, thanks to three decades of active research \cite{Yang2014},
and furthermore extreme sensitivity to the assumed model is not a
desirable property of an inference method.  Seven decades ago, Egon
S.\ Pearson \cite{Pearson1947} wrote that ``Hitherto the user has been
accustomed to accept the function of probability theory laid down by
the mathematicians; but it would be good if he could take a larger
share in formulating himself what are the practical requirements that
the theory should satisfy in application.'' This stipulation may be
relevant even today.

Two heuristic approaches have been suggested to remedy the high posterior model
probabilities in the context of phylogenies.  The first is to assign nonzero
probabilities to multifurcating trees (such as the star tree of Fig.~1) in
the prior \cite{Lewis2005}.  This is equivalent to assigning some prior probability to the
model $p = 0.5$ in the fair-coin example of Problem 1.  While this resolves the star-tree
parodox, it suffers from the conceptual difficulty that the multifurcating trees may not be
plausible biologically.  The second approach is to let the internal branch lengths
in the binary trees become increasingly smaller in the prior when the data size increases
\cite{Yang2005, Yang2007}. This is non-Bayesian in that the prior depends on the size of
the data.  With both approaches, the posterior is extremely sensitive to the prior \cite{Yang2008}.

\subsection{Non-Bayesian methods}

The phylogeny problem was described  by Jerzy Neyman \cite{Neyman1971} as ``a source of
novel statistical problems''.  In the Frequentist framework, test of phylogeny, or test of
nonnested models in general, offers challenging inference problems.  Note that in many
model selection problems, the model itself is not the focus of interest.  For example,
when an experiment is conducted to evaluate the effect of a new fertilizer, the
sensitivity of the inference to the assumed normal distribution with homogeneous variance
may be of concern, but the focus is not on the normal distribution itself.  In
phylogenetics, the phylogeny (which is a model) is of primary interest, far more important
than the branch lengths (which are parameters in the model).  Test of phylogeny is thus
more akin to significance/hypothesis testing than to model selection.  Model-selection
criteria such as AIC \cite{Akaike1973} or BIC \cite{Schwarz1978} simply rank the trees by
their likelihood (maximized over branch lengths), and will not be useful for attaching a
measure of significance or confidence in the estimated tree.  The phylogeny problem (or the
problem of comparing nonnested models in general) falls outside the Fisher-Neyman-Pearson
framework of hypothesis testing, which involves two nested models, one of which is true
\cite{Lehmann1997, Goldman2000}.

In principle Cox's likelihood ratio test \cite{Cox1961}, which conducts multiple tests
with each model used as the null, can be used to compare nonnested models.  For type-3
problems (Fig.~2, $C_1$-$C_4$), this test should lead to rejection of all
models.  Cox's test has not been used widely in phylogenetics, apparently because of  the
existence of a great many possible trees and the heavy computation needed to generate the
null distribution by simulation.

The most commonly used method for attaching a measure of confidence in the maximum
likelihood tree is the bootstrap \cite{Felsenstein1985}, which samples sites (alignment
columns) to generate bootstrap pseudo-datasets and calculates the bootstrap support value
for a clade (a node on the species tree) as the proportion of the pseudo-datasets in which
that node is found in the inferred ML tree.  This application of bootstrap for model
comparison appears to have important differences from the conventional bootstrap for
calculating the standard errors and confidence intervals for a parameter estimate
\cite{Efron1993}: a straightforward interpretation of the bootstrap support values for
trees remains elusive \cite{Felsenstein1993, Efron1996, Susko2009, Yang2014}.  At any
rate, the asymptotic behavior of bootstrap support values under the different scenarios of
Fig.~2 merits further research.  For the fair-coin example of problem 1
(Fig.~2, $C_1$), the bootstrap support converges to $U(0, 1)$, unlike the 
posterior probability, although other cases are yet to be explored.



\begin{materials}
Star-tree simulations.  For Fig.~4A, A', B \& B', the true tree is $T_0$ of Fig.~1A. The
data of counts of five site patterns ($xxx$, $xxy$, $yxx$, $xyx$, and $xyz$) were
simulated by multinomial sampling \cite{Yang1994SB}, and analyzed using a C program, which
calculates the 2-D integrals in the marginal likelihood by Gaussian-Legendre quaduature
with 128 points \cite{Yang2007}.  For Fig.~4C \& C', the true tree is $T_0$ of Fig.~1B. 
Sequence alignments were simulated using \textsc{evolver} and analyzed using MrBayes
\cite{Ronquist2012}.

\end{materials}




\begin{acknowledgments}
We thank Philip Dawid and Wally Gilks for stimulating discussions, and Jeff Thorne and an
anonymous reviewer for constructive comments.  Z.Y.\ was supported by a Biotechnological
and Biological Sciences Research Council grant (BB/P006493/1), and in part by the
Radcliffe Institute for Advanced Study at Harvard University.  T.Z.\ was supported by
Natural Science Foundation of China grants (31671370, 31301093, 11201224 and 11301294) and
a grant from the Youth Innovation  Promotion Association of Chinese Academy of Sciences
(2015080).
\end{acknowledgments}

\begin{table}
\caption{Proportions of datasets with extreme posterior probabilities for the
	three binary trees in the star-tree simulation.  }

\begin{tabular}{l c c c c c c}
\hline
\textbf{$n$} & $\mathbb{P}\{P_{\min}<1\%\}$ & $\mathbb{P}\{P_{\min}<5\%\}$ & $\mathbb{P}\{P_{\max}>95\%\}$ & $\mathbb{P}\{P_{\max}>99\%\}$ & $E(P_{\min})$ & $E(P_{\max})$ \\
\hline
$10^3$  & 0.234	& 0.550	& 0.205	& 0.079	& 0.067 & 0.754	 \\
$10^4$  & 0.812	& 0.931	& 0.606	& 0.450	& 0.011 & 0.897	 \\
$10^5$  & 0.979	& 0.992	& 0.853	& 0.773	& 0.001 & 0.964	 \\
$10^6$  & 1.000	& 1.000	& 0.953	& 0.919	& 0.000 & 0.988	 \\
$10^7$  & 0.999	& 1.000	& 0.982	& 0.970	& 0.000 & 0.995	 \\
\hline
\end{tabular}
$P_{\max} = \max(P_1, P_2, P_3)$ and
$P_{\min} = \min(P_1, P_2, P_3)$. Data are generated under JC+$\Gamma$ with $\alpha = 1$
using the star tree for four species: $(a: 0.2, b: 0.2, c: 0.2, d: 0.2)$, and analyzed
under JC.  The number of replicates is $10^3$. The probability density of $(P_1, P_2,
P_3)$ for the case of $n = 10^3$ and $10^5$ are shown in Fig.~4C and 4C'. 
\end{table}

\begin{table}
\caption{Proportions of datasets with strong support for wrong trees in
	simulated datasets for four species. }

	\begin{tabular}{l c c c c}
		\hline
		\textbf{$n$} & $\mathbb{P}\{P_1<1\%\}$ & $\mathbb{P}\{P_1<5\%\}$ & $\mathbb{P}\{P_{23}>0.95\%\}$ & $\mathbb{P}\{P_{23}>99\%\}$ \\
		\hline
		$10^3$  & 0.083	& 0.225	& 0.113	& 0.038 \\
		$10^4$  & 0.250	& 0.337	& 0.266	& 0.166 \\
		$10^5$  & 0.102	& 0.120	& 0.115	& 0.097 \\
		$10^6$  & 0.000	& 0.000	& 0.000	& 0.000 \\
		$10^7$  & 0.000	& 0.000	& 0.000	& 0.000 \\
		\hline
	\end{tabular}
$P_{23} = \max(P_2, P_3)$.  Data are generated under
JC+$\Gamma$ (with $\alpha = 1$) and analyzed under JC.  The true tree  is $T_1$ of
Fig.~1B: $((a:0.2, b:0.2):0.002, c:0.2, d:0.2)$, so that $T_2$ and $T_3$ are the two
wrong trees in the analysis.  The number of replicates is $10^3$. 
\end{table}





\end{article}

\end{document}


\begin{article}

\section{Supplemental Information Appendix}

\subsection{General theory for equally wrong models with no free parameters ($d = 0$)} 

The data, $\vec{x} = \{x_i\}$, consist of an i.i.d.\ sample from the true model
$g(\cdot)$; that is, $x_i \sim g(x_i)$, $i = 1, 2, \cdots, n$.  We consider two models
$H_1$ and $H_2$, with densities $f_1(x)$ and $f_2(x)$, respectively.  The models are
equally wrong, with $D_1 = D_2 > 0$ in eq.~[1] in the main text, and also identifiable.  We assign a
uniform prior ($\pi_k = \frac{1}{2}$  each) on the two models.  With no parameters in
either model, the marginal likelihood ($M$) is the same as the likelihood ($L$), so that
the logarithm of the marginal likelihood ratio is\\
\begin{equation}  \label{eq:A1}
   \Delta _n = \log \frac{M_1}{M_2} = \log \frac{L_1}{L_2} = \sum_i^n \log \frac{f_1(x_i)}{f_2(x_i)} = \sum_i^n r_i .
\end{equation}
Thus $\Delta_n$ is a sum of $n$ i.i.d.\ random variables ($r_i$).  

For large $n$, $\Delta_n$ has approximately a normal distribution by the central-limit theorem, with mean
\begin{equation} \label{eq:A2}
  \mathbb{E}_g(\Delta_n) = \mathbb{E}_g \left\{ \sum_i^n r_i \right\}
  = n \int g(x) \log \frac{f_1(x)}{f_2(x)} \textrm{d} x  = n(D_1 - D_2) = 0
\end{equation}
and variance
\begin{equation}  \label{eq:A3}
 \mathbb{V}_g(\Delta_n) = \mathbb{E}_g \left\{ \sum_i^n r_i^2 \right\}
  = n \int g(x) \left[ \log \frac{f_1(x)}{f_2(x)} \right] ^2 \textrm{d} x = nC,
\end{equation}
with $C > 0$, since the two models are distinct.  Note that $P_1 = \frac{1}{1 +
	1/\textrm{e}^{\Delta_n}}$. For $P_1$ to be not extreme, $\Delta_n$ should be close to 0. 
$P_1$ is in the interval $(\alpha, 1 - \alpha)$ for small $\alpha$, if and only if
$|\Delta _n| < A = \log \frac{1 - \alpha}{\alpha}$.  With large $n$, this occurs with probability
\begin{equation}  \label{eq:A4}
 \mathbb{P} \{ |\Delta_n| < A \} \approx 1 - 2\Phi \left( -\frac{A}{\sqrt {nC}} \right)  \approx \frac{2A}{\sqrt {2\pi nC}}.
\end{equation}

In Problem 1 (fair coin paradox), $r_i$ has a two-point distribution taking values $\pm
\log \frac{0.4}{0.6}$, each with probability $\frac{1}{2}$, so that $\Delta_n$ ($n = 0, 1,
\cdots$) constitutes a discrete-step random walk. We have $\mathbb{E}_g(\Delta_n) = 0$ and
$\mathbb{V}_g(\Delta_n) = nC = n \left[\log \frac{0.4}{0.6}\right]^2$, and eq.~\eqref{eq:A4}
agrees with eq.~3 in the main text.

We are interested in the Frequentist properties of Bayesian model selection.  If we
generate many replicate datasets under the true model $g(x)$ and analyze each to calculate
$P_1$, the proportion of datasets in which $P_1$ lies in the interval $(\alpha, 1 -
\alpha)$ goes to 0 at the rate of $n^{-\frac{1}{2}}$.  In the limit when $n \rightarrow
\infty$, $P_1 \rightarrow 0$ in half of the datasets and $\rightarrow 1$ in the other
half.  Previously Berk \cite{Berk1966} discussed the case of two equally wrong models
represented by $\theta = 0$ and 1, noting that asymptotically the posterior model probability
does not converge to a point value.

We say that $\Delta_n$ is of order $n^{\frac{1}{2}}$, or $\Delta_n =
\Theta_p(n^{\frac{1}{2}})$.  Formally $Y_n = \Theta_p(a_n)$ if for any given probability
$\epsilon > 0$, there exist $N$, $A_1(N, \epsilon) > 0$ and $A_2(N,\epsilon) > 0$, such
that when $n > N$, we have
\begin{equation}  \label{eq:A5}
\mathbb{P} \left\{ A_1 < |Y_n/a_n| < A_2 \right\} > 1 -\epsilon.
\end{equation}
Effectively $Y_n$ increases with $n$ no faster than $a_n$ and no more slowly than $a_n$. 
Using eqs.~\eqref{eq:A2} and \eqref{eq:A3}, with $Y_n = \Delta_n$ and $a_n =
\sqrt{n}$, it is easy to confirm that eq.~\eqref{eq:A5} is satisfied if $A_1 =
C^{\frac{1}{2}} \Phi^{-1}(\frac{1}{2} + \frac{\epsilon}{4})$ and $A_2 = C^{\frac{1}{2}}
\Phi^{-1}(1 - \frac{\epsilon}{4})$.  Thus when $n$ increases, $|\Delta_n|$ increases no faster
than $\sqrt{n}$ and no more slowly than $\sqrt{n}$.

In the case of comparing $K$ equally wrong models the above argument applies to any pair of models.  With probability $1/K$ (i.e., in $1/K$ of the datasets), the posterior model probability for one of the models 
$\rightarrow 1$ while all others $\rightarrow 0$, when $n \rightarrow \infty$.

Note the assumption that the two models are equally wrong, with $E_g(\Delta_n) = 0$ or
$D_1 = D_2 > 0$.  Otherwise if $D_1 \ne D_2$, $\Delta_n = n(D_1 - D_2)$ is $\Theta_p(n)$
and the less wrong model will dominate.

\subsection{General theory for equally right or equally wrong models with free parameters $(d>0)$}

The data ($\vec{x}$) are generated from the density $g(\cdot)$, and we compare two models
$H_1$ and $H_2$. Model $H_1$ specifies the density $f_1(x|\theta_1)$ with $d_1$ parameters
($\theta_1$), while $H_2$ has density $f_2(x|\theta_2)$ with $d_2$ parameters ($\theta_2$). 
We assign a uniform prior ($\pi_k = \frac{1}{2}$ each) on the two models, and the prior $f_k(\theta_k)$ for parameter $\theta_k$ under model $H_k$.
For any dataset, $x = \{x_1, \cdots, x_n\}$, the MLE $\hat{\theta}_k$  maximizes the
likelihood function $f_k(x|\theta_k)$ under model $H_k$.   When $n \rightarrow \infty$,
$\hat{\theta}_k \rightarrow \theta_k^*$. Thus $\hat{\theta}_k$  may be considered a
natural estimate of $\theta_k^*$  \cite{White1982}.  As usual, we assume that both
$\hat{\theta}_k$ and $\theta_k^*$ are inner points in the parameter space of $H_k$.  
Whether $\theta_k^*$ is inside the parameter space or at its boundary should affect
the precise distribution of $P_1$ but not its dynamics (i.e., whether or not $P_1$ has a
degenerate distribution).

As in \cite{White1982}, we define two matrices:
\begin{align} \label{eq:A6}
 {I_k}({\theta _k}) &= \mathbb{E}_g \{ \nabla \log f_k(x|\theta _k) \cdot \nabla \log f_k(x|\theta _k)^T \}, \nonumber \\
 {J_k}({\theta _k}) &= \mathbb{E}_g \{ - \nabla ^2 \log f_k(x|\theta_k) \},
\end{align}
where the expectation is over the true distribution $x \sim g(\cdot)$, and where $\nabla$
and $\nabla^2$  are the first and second derivatives with respect to $\theta_k$.

Following \cite{Dawid2011}, we decompose the marginal likelihood $M_k = f_k(x)$, $k = 1, 2$,
as a product of three terms, so that
\begin{equation}  \label{eq:A7}
 \log \frac{M_k}{g(x)} = \log \frac{f_k(x)}{f_k(x|\hat{\theta}_k)} 
                       + \log \frac{f_k(x|\hat{\theta}_k)}{f_k(x|\theta^*_k)} 
                       + \log \frac{f_k(x|\theta^*_k)} {g(x)}
                      := A_k + B_k + C_k.
\end{equation}
We define the corresponding differences between the two models as $\Delta A = A_1-A_2,
\Delta B = B_1-B_2$, and $\Delta C = C_1-C_2$, with
\begin{equation} \label{eq:A8}
   \Delta := \log \frac{M_1}{M_2} = \Delta A + \Delta B + \Delta C.
\end{equation}

From eq.~\eqref{eq:A6} of \cite{Dawid2011}, the first term in eq.~\eqref{eq:A8} is
\begin{equation} \label{eq:A9}
\Delta A = -\frac{d_1 - d_2}{2} \log \frac{n}{2\pi}
         + \log \frac{f_1(\theta_1^*)}{f_2(\theta_2^*)}
         + \log \left( \frac{\det J_2^*}{\det J_1^*} \right)^{\frac{1}{2}}
         + \Theta_p(n^{-\frac{1}{2}}),
\end{equation}
where $J_k^* = J_k(\theta_k^*)$ is $J_k(\theta_k)$ evaluated at $\theta_k^*$, and $\det Z$
is the determinant of matrix $Z$.

For the second term in eq.~\eqref{eq:A8}, $\Delta B$, we have
\begin{equation}  \label{eq:A10}
  B_k \approx \frac{1}{2} \left\{ \sqrt{n}(\hat{\theta}_k - \theta^*_k) \right\}^T   J^*_k 
                          \left\{ \sqrt{n}(\hat{\theta}_k - \theta^*_k) \right\}
\end{equation}
\cite[equation A.8 in]{Dawid2011}.  As $\sqrt{n}(\hat{\theta}_k - \theta_k^*)$ converges in
distribution to $\mathbb{N}\left(0, [(J_k^*)^{-1}]^T I_k^* (J_k^*)^{-1} \right)$
\cite[Theorem 3.2]{White1982}, $B_k$ and thus $\Delta B$ are $\Theta_p(1)$. Here $I_k^* =
{I_k}(\theta _k^*)$ is $I_k(\theta_k)$ evaluated at $\theta_k^*$ .  In the special case
that model $H_k$ is correct, $I_k^* = J_k^*$  and $B_k$ converges to $\frac{1}{2}
\chi_{d_k}^2$ in distribution.  $\Delta B$ is then the difference of two (correlated)
$\frac{1}{2} \chi_{d_k}^2$ variables.

Finally, the third term in eq.~\eqref{eq:A8}, $\Delta C = \log \frac{f_1(x|\theta_1^*)}
{f_2(x|\theta_2^*)}$ is identically 0 if the two models are indistinct, that is, if
${f_1}(x|\theta _1^*) = {f_2}(x|\theta _2^*)$ almost everywhere, as in the type-1 and type-2
problems of Fig.~2.  Otherwise $\Delta C$ is of
$\Theta_p(n^{\frac{1}{2}})$, given by the random walk (eq.~\eqref{eq:A1}).  As in the case
of no parameters discussed above, we assume that the two models are equally wrong, with
$D_1 = D_2$ or $\mathbb{E}_g(\Delta C) = 0$. Otherwise $\Delta C$ is $\Theta_p(n)$ and
dominates $\Delta A$ and $\Delta B$: the less wrong model will dominate with posterior
probability approaching 1.

The order of the three terms in eq.~\eqref{eq:A8} is summarized in table S1.  First is the
case where the two models are indistinct.  We have $\Delta C = 0$ and $\Delta B =
\Theta_p(1)$, while $\Delta A = \Theta_p(1)$ if $d_1 = d_2$ and $\Delta A =
\Theta(\log(n))$ if $d_1 \neq d_2$.  If the two models have the same dimension ($d_1 =
d_2$), $\Delta$ is of order $\Theta_p(1)$ and converges to a non-degenerate distribution,
which is determined by both $\Delta A$ and $\Delta B$.  This is the type-2 behavior
of Fig.~2.  If $d_1 \neq d_2$, the $\log n$ term in
eq.~\eqref{eq:A9} dominates, so that the model with fewer parameters wins. It is
noteworthy that as long as the two models are indistinct (and no matter whether they are
equally right or equally wrong), the model with fewer parameters wins.

Next is the case where the two models are distinct.  We have $\Delta A$ is $\Theta(1)$ if
$d_1 = d_2$ or $\Theta(\log(n))$ if $d_1 \ne d_2$, $\Delta B$ is $\Theta_p(1)$ and $\Delta
C$ is $\Theta_p(\sqrt{n})$, so that $\Delta C$ dominates.  Whether or not the two models
have the same dimension, $\Delta$ behaves like a random walk with mean 0 and variance of
order $n$.  In this case, $P_1 \rightarrow 1$ in half of the datasets and $\rightarrow 0$
in the other half.  The case with $d_1 = d_2$ is illustrated as the type-3 behavior of
Fig.~2.  Note that when the two models are equally wrong and distinct, the
size of the model does not matter and the model with fewer parameters does not dominate.

\subsection{Analysis of Problem 2 (two equally right models or equally wrong but indistinct models)}

The true model is $\mathbb{N}(0, 1)$, and we compare two models $H_1: \mathbb{N}(\mu,
1/\tau), \mu < 0$ and $H_2: \mathbb{N}(\mu, 1/\tau), \mu > 0$, with $\tau$ given.  The MLE
in infinite data is $\mu^* = 0$ in each model, and the two models are indistinct.  The
data ($\vec{x}$) are summarized as the sample mean $\bar{x}$.  We assign uniform prior
(with 1/2 each) for the two models, and $\mu \sim \mathbb{N}(0, 1/\xi)$, with $\xi$ given,
truncated to the appropriate range under each model.  The posterior model probability
$P_1$ can be derived by considering the posterior of $\mu$ under the model
$\mathbb{N}(\mu, 1/\tau)$ with $-\infty < \mu < \infty$.  As the prior precision is $\xi$
and the data (likelihood) precision is $n\tau$, the posterior of $\mu$ is $\mu|x \sim
\mathbb{N}\left( \frac{n\tau \bar x}{n\tau  + \xi}, \frac{1}{n\tau  + \xi} \right)$. Thus
\begin{equation}  \label{eq:A11}
  P_1 = \mathbb{P}\{H_1 | \vec{x}\} = \mathbb{P}\{ \mu <0| \vec{x} \} 
      = \Phi \left( - \frac{n\tau \bar x}{\sqrt{n\tau  + \xi}} \right).
\end{equation}
As $\bar{x}$ varies among datasets according to $\mathbb{N}(0, 1/n)$, 
we have $P_1 \approx \Phi \left(z \sqrt{\tau}\right)$ if $n$ is large, 
where $z =  -\sqrt{n} \bar{x} \sim \mathbb{N}(0, 1)$. 
The density of $P_1$ is given by a variable transform.  Note that 
$\left| \frac{\textrm{d} P_1} {\textrm{d} \bar x} \right| 
= \phi \left( -\frac{n\tau \bar x} {\sqrt{n\tau  + \xi}} \right) \times \frac{n\tau}{\sqrt{n\tau + \xi}} $, 
where $\phi(x)$ is the probability density function (PDF) for $\mathbb{N}(0, 1)$, and 
$\bar{x} = -\Phi^{-1}(P_1) \frac{\sqrt{n\tau + \xi}} {n\tau}$.

\begin{align}  \label{eq:A12}
f(P_1) &= \phi (\bar{x}; 0, \frac{1}{n}) \left/ \left| \frac{\textrm{d} P_1} {\textrm{d} \bar x} \right| \right. \nonumber\\
       &= \frac{1}{\sqrt{2\pi/n}} \textrm{e}^{-\frac{1}{2} n\bar{x}^2}
       \times \sqrt{2\pi} \cdot \textrm{e}^{\frac{1}{2} \frac{(n\tau \bar{x})^2}{n\tau+\xi}}
       \times \frac{\sqrt{n\tau+\xi}}{n\tau}       \nonumber\\
       &= \frac{\sqrt{\tau+\xi/n}}{\tau} 
       \cdot \exp\left\{ \frac{n}{2} \bar{x}^2 \left[ \frac{n\tau^2}{n\tau+\xi} - 1  \right] \right\}  \nonumber\\
       &= \frac{\sqrt{\tau+\xi/n}}{\tau} 
       \cdot \exp\left\{ \frac{n}{2} [\Phi^{-1}(P_1)]^2 \cdot \frac{n\tau+\xi}{(n\tau)^2}
                \left[  \frac{n\tau^2}{n\tau+\xi} - 1 \right] \right\} \nonumber\\
       &= \frac{\sqrt{\tau+\xi/n}}{\tau} 
       \cdot \exp\left\{ \frac{1}{2} [\Phi^{-1}(P_1)]^2 
                \left[ 1 - \frac{1}{\tau} - \frac{\xi}{n\tau^2} \right] \right\},
\end{align}
where $\phi(x; \mu, \sigma^2)$ is the PDF for $x \sim \mathbb{N}(\mu, \sigma^2)$.  This is eq.~4 in the main text.

\subsection{Analysis of Problem 3 (two equally wrong and distinct models, Gaussian with incorrect variances)}

Suppose the true model is $\mathbb{N}(0, 1)$, and two compared models are
$H_1$: $N(\mu, 1/\tau_1)$ and $H_2: \mathbb{N}(\mu, 1/\tau_2)$, with $\tau_1 < 1 <\tau_2$ given,
while $\mu$ is a free parameter.  The K-L divergence from model $H_1$ with parameter $\mu$ 
to the true model is
\begin{equation}  \label{eq:A13}
   D_1(\mu) = \int \phi(x) \log \frac{\phi(x)}{\phi(x; \mu, 1/\tau _1)} \textrm{d}x.
\end{equation}
$D_1(\mu)$ is minimized at $\mu^* = 0$.  Similarly for model $H_2$, $D_2(\mu)$ is
minimized at $\mu^* = 0$.  Letting $D_1(\mu^*) = D_2(\mu^*)$ so that the two models are
equally wrong leads to $\log \frac{\tau_1}{\tau_2} = (\tau_1 -\tau_2)$. In general, if the
true model is $\mathbb{N}(\mu, 1/\tau_0)$, then $H_1$ and $H_2$ are equally wrong if
$\tau_0 = (\tau_1 - \tau_2) \left/ \log \frac{\tau_1}{\tau_2} \right.$.

We assign a uniform prior (1/2 each) for the two models, and $\mu \sim \mathbb{N}(0,
1/\xi)$ under each model, with $\xi$ given.  The data are summarized as the sample mean
$\bar{x}$ and sample variance $ s^2 = \frac{1}{n} \sum_i (x_i - \bar{x})^2 $.  The
marginal likelihood under $H_1$ is
\begin{align}   \label{eq:A14}
  M_1 &= \int_{-\infty}^\infty \frac{1}{\sqrt{2\pi/\xi}} \textrm{e}^{-\frac{\xi}{2}\mu^2}
         \cdot \left(\frac{1}{\sqrt{2\pi/\tau_1}} \right)^n
         \exp\{ - \frac{1}{2}\tau_1 \sum_i (x_i - \mu )^2 \} \, \textrm{d}\mu  \nonumber \\
      &= \sqrt{ \frac{\xi}{2\pi} } \left( \frac{\tau_1}{2\pi} \right)^{\frac{n}{2}} 
         \int_{-\infty}^\infty \exp\left\{-\frac{1}{2}(\xi +n\tau_1)\mu^2 
                                  + n\tau_1\bar{x}\mu - \frac{1}{2}\tau_1\sum x_i^2 \right\} \textrm{d}\mu  \nonumber \\
      &= \sqrt{ \frac{\xi}{2\pi} } \left( \frac{\tau_1}{2\pi} \right)^{\frac{n}{2}} 
         \times \sqrt{ \frac{2\pi}{\xi+n\tau_1} }
         \exp\left\{  -\frac{1}{2}n\tau_1 
         \left( \bar{x}^2 + s^2 -\frac{n\tau_1\bar{x}^2}{\xi+n\tau_1} \right)
          \right\}                        
         \nonumber \\
      &= \sqrt{ \frac{\xi}{\xi+n\tau_1} } \left( \frac{\tau_1}{2\pi} \right)^{\frac{n}{2}} 
       \times \exp\left\{ -\frac{n\tau_1}{2(\xi+n\tau_1)} [\xi\bar{x}^2 + \xi s^2 + n\tau_1 s^2] 
       \right\}.
\end{align}

With the marginal likelihood $M_2$ under $H_2$ similarly given, the posterior odds is
\begin{align} \label{eq:A15}
 & \frac{P_1}{P_2} = \sqrt{ \frac{\xi+n\tau_2}{\xi+n\tau_1} }
                       \left[ \frac{\tau_1}{\tau_2} \right] ^ {\frac{n}{2}} 
                       \exp \left\{\frac{n}{2(\xi+n\tau_1)(\xi+n\tau_2)}   \right. \nonumber \\
  &   \times \left. \left[ (\tau_2-\tau_1)(\xi^2\bar x^2 + \xi^2 s^2 + n^2\tau_1\tau_2 s^2) + (\tau_2^2-\tau_1^2)n\xi s^2 \right] \right\}.
\end{align}

The distribution of $P_1$ is given as a transform of $\bar x$ and $s^2$, but the
derivation is tedious.  Instead we generate a sample of $P_1$ by simulating random
variables $\bar x \sim \mathbb{N}(0, 1/n)$ and $ns^2 \sim \chi_{n-1}^2$.

When $n \rightarrow \infty$, we obtain, using eq.~5 in the main text,
\begin{equation}  \label{eq:A16}
  \frac{P_1}{P_2} = \frac{M_1}{M_2} \rightarrow \textrm{e}^{ \frac{1}{2}(\tau_2 - \tau_1)(n s^2 - (n-1)) }
                  \rightarrow \textrm{e}^{ \frac{1}{2}(\tau_2-\tau_1) Z } ,
\end{equation}
where $Z = ns^2 - (n-1) \sim \mathbb{N}(0, 2n-2)$ because for large $n$, $\chi_{n-1}^2$ is
approximately $\mathbb{N}(n-1, 2n-2)$.  Thus when $n$ is large, $Z$ will have a
vanishingly small probability of lying in a fixed interval around 0, and $P_1$ will be
either $\sim 0$ or $\sim 1$, each in half of the datasets.

\subsection{Star-tree simulation} 

We consider three cases, involving trees of three or four species only.  

In case 1 (Fig.~4A \&A'), the molecular clock (rate constancy over time)
is assumed, and the true tree is the star tree for three species of Fig.~1A,
with the branch length $t = 0.2$ (changes per site on average).  The JC substitution model
\cite{Jukes1969} is used both to generate and to analyze the data.  The three binary trees $(T_1, T_2,
T_3)$ represent three equally right models.  We assign a uniform prior ($\frac{1}{3}$
each) on the binary trees, and given each binary tree, exponential priors with means 0.1
for $t_0$ and 0.2 for $t_1$.  The marginal likelihoods or the 2-D integrals over $t_0$ and
$t_1$ are calculated numerically using Gaussian quaduature \cite{Yang2007}.

Case 2 (Fig.~4B \&B') is similar to case 1, except that the substitution model used in the
simulation is JC+$\Gamma$ \cite{Jukes1969, Yang1993}, with different sites in the sequence
evolving at variable rates according to the gamma distribution with shape parameter
$\alpha = 1$.  The analysis model is still JC (equivalently to JC+$\Gamma$ with $\alpha =
\infty$). The MLE at the limit of infinite data is $t_0^*  = 0$ in each binary tree, so
that the three models under comparison are equally wrong and indistinct.

In case 3 (Fig.~4C \&C') the true model is JC+$\Gamma$ with $\alpha = 1$ and the analysis
model is JC, as in case 2, but we simulate data using an unrooted star tree for four
species.  The true branch lengths are $t_1 = t_2 = t_3 = t_4 = 0.2$ in $T_0$ (Fig.~1B). 
The data are analyzed to evaluate the three binary unrooted trees for four species
(Fig.~1B), each with five branch lengths, without assuming the molecular clock.  Because
the MLE $t_0^{*} > 0$  in the binary trees in the limit of infinite data, this case is
a comparison of equally wrong and distinct models.  The three binary trees are assigned
equal prior probabilities ($\frac{1}{3}$ each), and the five branch lengths in each tree
are assigned the exponential prior with mean 0.1.  The program MrBayes \cite{Ronquist2012}
is used to calculate the posterior probabilities for the three binary trees $(P_1, P_2,
P_3)$ by MCMC.  The program is tricked into collapsing the $n$ sites into 15 (instead of
256) site patterns to speed up the likelihood calculation.  Each MCMC run, for $2 \times
10^6$ iterations, takes $\sim 10$ seconds.

\begin{table*} [h]
	\caption{Table S1.  The order of the three terms ($\Delta A, \Delta B$, and $\Delta C$) in eq.~\eqref{eq:A8} and the asymptotic behavior of Bayesian model selection}  \label{tableconvergencerate}
	\begin{tabular}{ l c c c l }
		\hline
& $\Delta A$ & $\Delta B$ & $\Delta C$ & Behavior of Bayesian model selection \\
\hline \\
\multicolumn{5}{l}{Indistinct models } \\
$d_1 = d_2$ & $\Theta_p$(1) & $\Theta_p(1)$  & 0 & Converges to a nondegenerate distribution \\
$d_1 \ne d_2$ & $\Theta(\log(n))$ & $\Theta_p(1)$  & 0 & Model with fewer parameters dominates \\
\\
\multicolumn{5}{l}{Distinct models} \\
$d_1 = d_2$  & $\Theta_p(1)$& $\Theta_p(1)$  & $\Theta_p(n^{\frac{1}{2}})$ & Random walk \\
$d_1 \ne d_2$ & $\Theta_p(\log{n})$ & $\Theta_p(1)$ & $\Theta_p(n^{\frac{1}{2}})$ & Random walk \\
\hline
\end{tabular}
\end{table*}

\end{article}